\documentclass[11pt]{amsart}
\usepackage{amscd,amssymb}

\numberwithin{equation}{section}

\newtheorem {Theorem}[equation]     	{Theorem}
\newtheorem {Lemma}[equation]     	{Lemma}

\newtheorem {proposition}[equation]	{Proposition} 
\newtheorem {Proposition}[equation]	{Proposition}

\newtheorem*{th1}{Theorem~1}
\newtheorem*{th1'}{Theorem~$\mathbf 1'$}

\theoremstyle{definition}

\newtheorem {definition}[equation]{Definition}
\newtheorem {Remark}[equation]		{Remark}

\def	\inv	{^{-1}}

\def	\index	{\operatorname{index}}
\def	\Res	{\operatorname{Res}}

\def	\id	{\operatorname{id}}
\def    \red    {M_{\text{red}} }
\def \conec     {\overset \circ  c }

\newcommand{\fg} {{\mathfrak g}}

\def	\cF	{{\mathcal F}}

\def	\cU	{{\mathcal U}}

\def	\C	{{\mathbb C}}
\def	\bF	{{\mathbb F}}
\def	\R	{{\mathbb R}}
\def	\bN	{{\mathbb N}}

\def	\Z	{{\mathbb Z}}

\def	\tX	{{\tilde{X}}}

\def	\tred	{{ \widetilde{M}_{\text{red}} }}
\def	\tPhi	{{\tilde{\Phi}}}
\def    \tT     {\tilde{T} } 
\def    \trho     {\tilde{\rho} }

\def	\bp	{{\overline{p}}}

\def	\bm	{{\overline{m}}}

\begin{document}

\title{Intersection cohomology of $S^1$ symplectic quotients and small
resolutions}

\author{Eugene Lerman}\author{Susan Tolman}
\date{\today}

\address{Department of
Mathematics, University of Illinois, Urbana, IL 61801}
\email{lerman@math.uiuc.edu}\email{stolman@math.uiuc.edu} 

\begin{abstract}
We provided two explicit formulas for the intersection cohomology (as a
graded vector space with pairing) of the symplectic quotient by a
circle in terms of the $S^1$ equivariant cohomology of the original
symplectic manifold and the fixed point data.  The key idea is the
construction of a small resolution of the symplectic quotient.
\end{abstract}

\maketitle 

\section{Introduction}

Let a compact Lie group act effectively on a compact connected
symplectic manifold $M$ with a moment map $\Phi: M \to \fg^*$.  In the
case that $0$ is a regular value of the moment map, the symplectic
quotient $\red := \Phi \inv (0)/G$ is an orbifold, and its rational
cohomology ring is fairly well understood \cite{Kibook, Wi, Wu, Ka,
Ka-Gu, JK, TW}.

However, many interesting spaces arise as reduced spaces at singular
values of the moment map.  Some examples include: the moduli space of
flat connections, some polygon spaces, many physical systems, and
projective toric varieties whose polytopes are not simple.  Since the
symplectic quotient at a singular value is a stratified space
\cite{SjL}, a natural invariant to compute is the intersection
cohomology (with middle perversity).  Less is known in this case.
Kirwan has provided formulas to compute the Betti numbers in the
algebraic case \cite{Kibook, Ki4}; Woolf extended this work to the
symplectic case.  Moreover, Jeffrey and Kirwan computed the
pairing in the intersection cohomology of particular symplectic
quotients \cite{ktalk}.

The main result of this paper is two explicit formulas for the
intersection cohomology (as a graded vector space with pairing) of the
symplectic quotient by a circle in terms of the $S^1$ equivariant
cohomology of the original symplectic manifold and the fixed point
data.  More precisely, these formulas depend on the image of the
restriction map in equivariant cohomology $H^*_{S^1}(M; \R) \to
H^*_{S^1}(M^{S^1}; \R)$ from the original manifold to the fixed point
set.
  
\begin{th1}
Let the circle $S^1$ act on a compact connected symplectic manifold
$M$ with moment map $\Phi: M \to \R$ so that $0$ is in the interior of
$\Phi(M)$.  Let $\red := \Phi\inv(0)/S^1$ denote the reduced space.

Then there exists a surjective map $\kappa$ from the equivariant
cohomology ring $H^*_{S^1}(M;\R)$ to the intersection cohomology
$IH^*(\red ;\R)$.  Moreover, given any equivariant cohomology class
$\alpha$ and $\beta$ in $H^*_{S^1}(M)$, the pairing of
$\kappa(\alpha)$ and $\kappa(\beta)$ in $IH^*(\red)$ is given by the
formula
$$
\left< \kappa(\alpha), \kappa(\beta) \right>  = 
 \Res _0 \sum_{F \in \cF^+} \int_F \frac{ i_F^*(\alpha \beta) }{e_F} .
$$
Here, $e_F$ denotes the equivariant Euler class of the normal bundle
of $F$, and $\cF^+ $ denotes the set of components $F$ of the fixed
point set ${S^1}$ such that either
\begin{enumerate}
\item $\Phi(F) > 0$ or
\item  $\Phi(F) = 0$ and $ \index  F \leq \frac{1}{2} (\dim M -  \dim F)$,
\end{enumerate}
where the index of $F$ is the dimension of the negative eigenspace
of the Hessian of the moment map $\Phi$ at a point of
$F$.
\end{th1}

The meaning of the right hand side is as follows.  The map $i_F^*$ is
simply the restriction to $F$.  The equivariant cohomology ring
$H^*_{S^1}(F)$ is naturally isomorphic to the polynomial ring in one
variable $H^*(F)[t]$.  The equivariant Euler class $e_F$ is invertible
in the localized ring $H^*(F)(t)$; thus, $\frac{i_F^*(\alpha)}{e_F}$
is an element of this ring.  The integral $\int_F : H^*(F)(t) \to
\R(t)$ acts by integrating each coefficient in the series. Finally,
$\Res_0$ denotes the operator which returns the coefficient of
$t\inv$.

\begin{Remark}
Our convention is that the pairing in intersection cohomology
between two classes $\alpha \in IH^p(\red)$ and $\beta \in IH^q(\red)$
is zero if $p + q \neq \dim(\red)$.
\end{Remark}

Note that since $\kappa$ is surjective, this theorem determines the
pairing for all pairs of elements in $IH^*(\red)$.  Additionally, by
Poincare duality in intersection cohomology, it determines the kernel
of $\kappa$.

We now provide an alternative  version of our main result:

\begin{th1'}
Let the circle $S^1$ act on a compact connected symplectic manifold
$M$ with moment map $\Phi: M \to \R$ so that
$0$ is in the interior of $\Phi(M)$.
Let $\red := \Phi\inv(0)/S^1$
denote the reduced space.

Then there exists a ring structure on the intersection
cohomology $IH^*(\red;\R)$ so that
\begin{itemize}
\item The ring structure on $IH^*(\red)$ is compatible with the pairing,
in the sense that their exists an isomorphism $\int$ from the top
dimensional intersection cohomology to $\R$ so that
$\int \alpha \cdot \beta = \left< \alpha, \beta \right>$.
\item As a graded ring,  $IH^*(M;\R)$ is isomorphic to 
$H^*_{S^1}(M;\R)/K$, where 
$$ K:=
 \{ \alpha \in H_{S^1}^*(M) \mid \alpha|_{F} = 0 \ \ \forall \ \ 
F \in {\cF}^+\} \oplus
 \{ \alpha \in H_{S^1}^*(M) \mid \alpha|_{F} = 0 \ \ \forall \ \ 
F \in {\cF}^-\} .$$
\end{itemize}
 
Here, $\cF^+$ denotes the set of components  $F$
of the fixed point set $M^{S^1}$  such that either
\begin{enumerate}
\item $\Phi(F) > 0$ or
\item  $\Phi(F) = 0$ and $ \index  F \leq \frac{1}{2} (\dim M -  \dim F)$,
\end{enumerate}
where the index of $F$ is the dimension  of the negative eigenspace
of the Hessian of the moment map $\Phi$ at a point of
$F$.  Additionally, $\cF^-$ denotes the set of all
other components of the fixed point set. 
\end{th1'}

In principle,
these  two formulas for the intersection 
cohomology  give almost exactly  the same information.  
We include both, because, in
practice, one or the other might be better suited to tackle a
particular problem.

We prove these two theorems simultaneously.  First, we construct an
orbifold $\tred$, which we call the {\bf perturbed quotient}, The
perturbed quotient is a small resolution of the symplectic quotient;
thus. as a graded vector space with pairing, $IH^*(\Phi\inv (0)/S^1)$
is isomorphic to $H^*(\widetilde{\red})$. Moreover, even though the
perturbed quotient is not symplectic, it is constructed in such a way
that the standard techniques for computing the cohomology ring of a
symplectic quotient can be applied to it, yielding the above formulae.

The construction of the perturbed quotient is fairly straightforward.
 The singularities of the reduced space $\red := \Phi\inv(0)/S^1$
 correspond to components $Y$ of the fixed point set $M^{S^1}$ lying
 on the zero level set $\Phi\inv(0)$.  In the setting projective
 varieties, it is known that the neighborhoods of these singularities
 have small resolutions \cite{Hu}.\footnote{
It is claimed in \cite{Hu} that the resolutions are global,
which need not be the case; see \cite{Hu2},} 
Although these resolutions are only local,
it is possible to piece them together into a global resolution.

We construct the perturbed quotient as the quotient of a fiber
of a perturbation $\tPhi: M \to \R$ of the original moment map.  Since
this perturbed moment map $\tPhi$ is Bott-Morse, and since its
critical points are exactly the fixed points of the action of $S^1$ on
$M$, the standard techniques used to compute the
cohomology of symplectic quotients can also be applied to compute the 
cohomology of the perturbed quotient.

Finally, we construct a pairing preserving isomorphism between
the intersection cohohomology of the symplectic quotient and
the cohomology of the perturbed quotient.
In the algebraic case, this follows immediately from the
fact that the perturbed quotient is a small resolution of the
symplectic quotient (see \S 6.2 in \cite{GM2}).  However, while
it ``appears to be clear that this theorem will also be valid in
our case", we have decided to provide a direct proof.

\subsection*{\sc Acknowledgments}
The work in this paper was inspired by  the lectures of Francis Kirwan
at the Newton institute in the Fall of 1994.

We thank Reyer Sjamaar for many helpful discussion.  In particular the
idea that for $S^1$ quotients the intersection cohomology should be
very simple  to compute is due to him.  We thank Sam Evens for a
number of useful discussions.
\\

\section{Simple stratified spaces and  intersection cohomology} 

In this section, we  introduce the two main concepts that we
will need in this paper: simple stratified spaces and intersection
cohomology.  The notion of a simple stratified space is not standard;
it is, however,  convenient for our purposes.  The definition of
intersection cohomology we use is essentially identical to the
definition of the complex of intersection differential forms due to
Goresky and MacPherson (see \cite{B}), except that we allow the strata
to be orbifolds, and that we only consider simple stratified spaces.

Recall that an open {\bf cone} on a topological space $L$ is
$$
\overset \circ  c (L) := L\times [0, 1)/\sim,
$$ 
where $(x,0) \sim (x',0)$ for all $x, x'\in L$.  Equivalently $\conec
(L) = L\times [0, \infty)/\sim$.

\begin{definition}\label{def_sss} 
A {\bf simple stratified space} is a Hausdorff topological space $X$
with the following properties:
\begin{itemize}
\item The space $X$ is a disjoint (set-theoretic) union of  orbifolds,
called {\bf strata}.
\item There exists an open dense oriented stratum $X^r$, called the {\bf top
stratum}. 
\item The complement of $X^r$ in $X$ is a disjoint union of 
connected orbifolds, $X \smallsetminus X^r = \coprod Y_i$,
called the {\bf singular strata}.
\item For each singular stratum $Y$ there is a neighborhood $\tT$ of
$Y$ in $X$ and a map $\pi :\tT \to Y$ which is a $C^0$ fiber bundle
with a typical fiber $\overset \circ c (L)$ for some orbifold $L$,
which depends on $Y$.  (Thus $Y$ embeds into $\tT$ as the vertex
section.)
\item Their exists a diffeomorphism from the
complement $\tT \smallsetminus Y$ 
to $Q \times (0,1)$, where $ Q\to Y$ is a
$C^\infty$ fiber bundle with typical fiber $L$,  such that the following
diagram commutes:
$$
\begin{CD}
\tT \smallsetminus Y @>>> Q_i \times (0,1)\\
@V{\pi}VV	                    @VVV\\
Y  @= 			Y.
\end{CD}
$$
In particular $\pi :\tT \smallsetminus Y \to Y$ is a smooth
fiber bundle with a typical fiber $L \times (0,1).$
\end{itemize}
\end{definition}
Thus a simple stratified space $X$ is a decomposition $X= X^r \sqcup
\coprod Y_i$ and a collection of maps $\{\pi_i : \tT_i \to Y_i \}$.  

\begin{Remark}
	\label{rmrk_t}
Note that the composite $\tT\smallsetminus Y \to Q \times
(0,1)\to (0,1)$, where $Q\times (0,1)\to (0,1)$ is the obvious
projection, extends to a continuous map $r: \tT\to [0, 1)$.  In the
definition of intersection cohomology of $X$ it
will be convenient for us to consider smaller tubular neighborhoods
$T$ defined by 
$$
	T = r\inv ([0, 1/2)).
$$
\end{Remark}
\begin{definition}[Cartan filtration of forms relative to a
submersion] Let $\pi : E\to B$ be a smooth submersion of orbifolds.
The Cartan filtration ${\bF}_k\Omega^*(E)$ of the complex of forms
$\Omega^* (E)$ on $E$ is given by
\begin{equation*}
\begin{split}
{\bF}_k\Omega^*(E)=\{ \omega \in \Omega^* (E) \mid & \ \ \text{for all
}\, e\in E \, \text{and  for all
vectors } \xi_0, \cdots, \xi_{k} \in \ker d\pi _e \\
& \ \ i(\xi_0) \circ \cdots \circ i(\xi_{k})(\omega(e))  = 0 \\
& \ \ \text{and } i(\xi_0) \circ \cdots \circ i(\xi_{k}) (d\omega (e)) = 0  \} .
\end{split}
\end{equation*}
\end{definition}
\noindent
By convention, $i(\xi_0) \circ \cdots \circ i(\xi_{k})(\sigma)
= 0$ if $\deg \sigma \leq k$. Note that ${\bF}_0 \Omega ^*
(E) $ consists of basic forms.

Let $X= X^r \sqcup \coprod Y_i$ be a simple stratified space. A {\bf
perversity} $\bar{p} :\{Y_i\} \to \bN$ is a function that assigns a
nonnegative integer to each singular stratum $Y_i$.  The {\bf middle
perversity} $\bm$ is defined by
$$
\bar{m} (Y_i)= \lfloor \frac{1}{2} (\dim X^r -\dim Y_i) \rfloor  -1.
$$

\begin{definition} Let $(X= X^r \sqcup \coprod Y_i,\{\pi_i : \tT_i \to
Y_i \})$ be a simple stratified space and let $\bp : \{Y_i\} \to \bN$
be a perversity.  The {\bf complex of intersection differential forms}
$I\Omega_{\bar p}^* (X)$ is a sub-complex of the complex of
differential forms on the top stratum:
$$
I\Omega_{\bar p}^* (X):= \{ \omega \in \Omega^* ( X^r) \mid \omega
|_{T_i \cap X^r} \in {\bF}_{\bar{p} (Y_i)} \Omega^*(T_i \cap
X^r)\}
$$
where the filtration ${\bF}_{\bar{p} (Y_i)} \Omega^*(T_i \cap
X^r)$ is defined relative to the submersion $\pi_i: T_i \cap X^r = T_i
\smallsetminus Y_i \to Y_i$.  The coboundary map is the exterior
differentiation $d$.

The {\bf intersection cohomology} $IH_{\bar{p}}^* (X)$ of the simple
stratified space $X$ with perversity $\bar{p}$ is the cohomology of
the complex $(I\Omega _{\bar{p}}^* (X), d)$.
\end{definition}

\begin{Remark}
If the strata of a simple stratified space $X$ are manifolds, then $X$
is a pseudomanifold. In this case our definition of intersection forms
is exactly the Goresky-MacPherson definition of the complex of
differential forms and our definition of intersection cohomology
agrees with the standard definition (see \S~1.2 in \cite{B}).
\end{Remark}

We now define the pairing on the middle perversity intersection
cohomology of a compact simple stratified space $X$ with an oriented
top stratum. 

Note that if $q > \dim(Y_i) + \bp(Y_i)$, then
every $\alpha \in  I\Omega^q_\bp(X)$ vanishes on $T_i$.
In particular, if $\dim X^r >  \dim Y_i + \bp(Y_i)$ 
for all singular strata $Y_i$ then every 
$\alpha \in I\Omega_\bp^{\dim X^r} (X)$ is supported on the compact set
$X \smallsetminus \bigcup T_i$. 
Therefore, if the top stratum $X^r$ is  oriented,
there is a well-defined integration map $\int : I\Omega_\bp^{\dim X^r} (X) \to
\R$, $\alpha \mapsto \int _{X^r} \alpha$. 
Similarly, if $\dim X^r - 1 >  \dim Y_i + \bp(Y_i)$ for all $i$,
then any $\beta \in I\Omega_\bp^{\dim X^r - 1}(X)$
is also supported in $X \smallsetminus \bigcup T_i$.
Thus, integration descends to a
well-defined map on cohomology $\int : IH_\bp^{\dim X^r} (X) \to \R$;
we extend this by zero to a map $\int : IH_\bp^*(X) \to \R$.

Given $\alpha$ and $\beta$ in $I\Omega^*_\bm(X)$, 
notice that $\alpha \wedge \beta \in I\Omega_{2 \bm}(X)$, where $2
\bm$ is twice the middle perversity.  This follows from the property
of the Cartan filtration: for any $\alpha \in {\bF}_k\Omega^*(E)$ and
$\beta \in {\bF}_l\Omega^*(E)$, $\alpha \wedge \beta \in
{\bF}_{k+l}\Omega^*(E)$.  Moreover, $2 \bm(Y_i) \leq \dim X^r - \dim
Y_i - 2$ for all singular strata $Y_i$.  Thus,
there is a well-defined bilinear pairing $IH_\bm^p (X) \times IH_\bm^q
(X) \to \R$ which sends $[\alpha] \in IH_\bm^p (X)$ and $[\beta] \in
IH_\bm^q (X)$ to the integral $\int_{X^r} \alpha \wedge \beta$.

\section{The structure of the symplectic quotient}

In this section, we recall a normal form for the neighborhoods of
fixed points on symplectic manifolds with Hamiltonian circle actions.
Using this, we give a normal form for the neighborhoods of the
singularities in a symplectic quotient. In particular, we show that the
quotient is a simple stratified space.

This last statement is a special case of a theorem of Sjaamar and Lerman
 \cite{SjL},
who show that every symplectic quotient by a compact Lie group is
a stratified space.  Note, however, that in \cite{SjL} the  stratification 
is by orbit type, whereas here we use the slightly coarser
stratification by infinitesimal orbit type.

Let a circle act on a symplectic manifold $M$ in a Hamiltonian fashion
with a moment map $\Phi: M \to \R$.  Recall that the {\bf symplectic
quotient} (a.k.a.\ the {\bf reduced space}) is
$M_{\text{red}}:= \Phi\inv(0)/S^1$.  If $0$ is a regular value for
$\Phi$, then the quotient is a symplectic orbifold.
More generally, $\Phi$ is regular on $M\smallsetminus M^{S^1}$, and
$\red^r := \left(\mu\inv (0) \cap (M\smallsetminus M^{S^1})\right)/S^1$
is an orbifold; this is the top stratum.
Moreover, since  the restriction of the
symplectic form on $M$ to $\Phi\inv (0)\cap (M\smallsetminus M^{S^1})$
descends to a symplectic form on $\red^r$, $\red^r$ is naturally
oriented.  

Recall that the moment map is constant on each component of the fixed
point set $M^{S^1}$, and that these components are isolated.  Thus,
every component $Y$ of the fixed point set which intersects the zero
level set is entirely contained in the level set, and gives rise to a
stratum of $\red$ diffeomorphic to $Y$.  To see how these strata fit
together, we need the following lemma.

\begin{Lemma}\label{lemma_model}
Let $S^1$ act on a symplectic manifold $(M,\omega)$ in a Hamiltonian
fashion with a moment map $\Phi: M \to \R$.  Every connected
component $Y$ of the fixed point set $M^{S^1}$ has even
codimension, say $2n$.   Moreover, there exists
\begin{itemize}
\item  positive integers $n_1,\ldots,n_k$ such
that $\sum_i n_i =  n $,
\item a principal $G := \prod U(n_i) \subset U(n)$ bundle $P$ over $Y$, and
\item  distinct non-zero weights  $\kappa_1,\ldots,\kappa_k$, 
\end{itemize}
such that:
\begin{itemize}
\item there is a diffeomorphism $\sigma$ from a neighborhood $U$ of
$Y$ in $M$ to a neighborhood $U_0$ of the zero section in the associated
bundle $ P \times _G \C^n \to Y$;
\item
this diffeomorphism is equivariant with respect to the 
circle action on $P \times_G \C^n$ defined by the weights $\kappa_i$;
\item
the diffeomorphism  pulls back the moment map $\Phi$ to the map
$\mu: P \times_G \C^n \to \R$ given below, 
i.e., $\Phi \circ \sigma = \mu,$ where
$$
\mu ([q, (\vec{z}_1, \ldots, \vec{z}_k)]) = \frac{1}{2} \sum \kappa _i
|\vec{z}_i|^2 + \Phi (Y), \quad \forall \ (\vec{z}_1, \ldots,
\vec{z}_k)\in \C^{n_1} \oplus \cdots \oplus \C^{n_k} = \C^n.
$$
\end{itemize}
\end{Lemma}

\begin{proof}
Consider the symplectic perpendicular bundle $E = TY^\omega$ of $Y$ in
$(M, \omega)$.  Since $Y$ is a symplectic submanifold of $M$ we have
$TM|_Y = TY \oplus TY^\omega$. So $E$ is the normal bundle of $Y$, and
$E$ is a symplectic vector bundle.  The group $S^1$ acts on the bundle
$E$ by fiber-preserving vector bundle maps.  We may choose an $S^1$
invariant complex structure on $E$ compatible with the symplectic
structure.  Up to an equivariant homotopy, this complex structure is
unique.
 
A fiber $\C^n$ of $E$ splits into the direct sum of isotypical
representations of $S^1$, $\C^n = \C^{n_1}\oplus \cdots \C^{n_k}$, so
that the action of $\lambda \in S^1$ on $\C^{n_i}$ is given by
multiplication by $\lambda ^{\kappa_i}$ for some weight $\kappa_i \in
\Z$.  Under the above identification of the fiber of $E$ with $\C^n$
the symplectic structure is the imaginary part of the standard
Hermitian inner product.  Hence a moment map for the $S^1$ action on
the fiber is
$$
 \C^{n_1}\oplus \cdots \oplus \C^{n_k}\ni (\vec{z}_1, \ldots,
 \vec{z}_k) \mapsto \frac{1}{2} \sum \kappa _i |\vec{z}_i|^2.
$$

The structure group of the vector bundle $E$ reduces to the subgroup
of $U(n)$ consisting of  transformations that commute with the
action of $S^1$ described above, that is, it reduces to $G := \prod
U(n_i) \subset U(n)$.  Consequently $E = P\times _G \C^n$ for some
principal $G$ bundle $P$ over $Y$.

The equivariant symplectic embedding theorem (see, for example,
Theorem~2.2.1 in \cite{GLS} and the subsequent discussion) implies
that we may identify the neighborhood of the submanifold $Y$ in $M$
with a neighborhood of the zero section of $E$ in such a way that a
moment map $\mu : E \to \R$ is given by
$$
\mu ([p, (\vec{z}_1, \ldots, \vec{z}_k)]) = 
		\frac{1}{2} \sum \kappa _i |\vec{z}_i|^2 + \text{a constant}
$$
for all $([p, (\vec{z}_1, \ldots, \vec{z}_k)] \in P\times _G (\C^{n_1}\oplus
\cdots \oplus \C^{n_k}).$ 
\end{proof}

\begin{Remark}\label{split} 
Let  $V^+$ and $V^-$ be the the sum of the positive and negative
weight spaces, respectively.
That is, we may assume that the corresponding weights satisfy 
$\kappa_1, \ldots \kappa_s >0$ and $\kappa _{s+1}, \ldots, \kappa _k < 0$, 
and set $$
V^+ = \bigoplus _{i=1}^s \C^{n_i}, \qquad 
V^- = \bigoplus _{i=s+1}^k \C^{n_i}.
$$
We then have linear representations of $G$ on $V^+$ and $V^-$
so that $E$ splits: $E = E^+ \oplus E^-$ where $E^\pm = P \times _G
V^\pm$.
By Lemma~\ref{lemma_model}, 
the index of the
moment map $\Phi :M\to \R$ at $Y$ is $\dim V^-$.
\end{Remark}

We now use Lemma~\ref{lemma_model} above to show that the reduced
space $\red$ is a simple stratified space.

\begin{proposition} \label{prop_simple}
Let the circle $S^1$ act effectively on a compact connected symplectic
manifold $M$ in a Hamiltonian fashion with a moment map $\Phi: M \to
\R$. Assume that 0 is in the interior of the image of the moment
map. The reduced space $\red = \Phi\inv (0)/S^1$ is a simple
stratified space.  The singular strata of $\red$ are
connected components $Y$ of the fixed point set $M^{S^1}$ with $\Phi
(Y) = 0$.  For each such stratum 
there exists:
\begin{itemize}
\item a faithful unitary representation $\rho: S^1 \to U(p) \times U(q)$,
where $p$ and $q$ are  positive integers whose sum is the codimension of $Y$;
and
\item a principal $G$ bundle $P \to  Y$, where $G$ 
is a subgroup of $U(p) \times U(q)$ which commutes with $\rho(S^1)$;
\end{itemize}
so that a neighborhood $\tilde{T}$ of $Y$ in $M$ is the associated cone bundle
$P \times_G \conec (S^{2p-1} \times_{S^1} S^{2q-1})$.
\end{proposition}
\begin{Remark}
By ``a neighborhood $\tilde{T}$ of $Y$ in $\red$ is the associated cone
bundle $P \times_G \conec (S^{2p-1} \times_{S^1} S^{2q-1})$'' we mean
that there exists a stratum-preserving homeomorphism from a
neighborhood $\tilde{T}$ of $Y$ in $\red$ to the associated bundle $P
\times_G \conec (S^{2p-1} \times_{S^1} S^{2q-1})$, which restricts to a
diffeomorphism on each stratum.  The stratification of $P \times_G
\conec (S^{2p-1} \times_{S^1} S^{2q-1})$ comes from the stratification
of the cone $\conec (S^{2p-1} \times_{S^1} S^{2q-1})$ into the vertex
and the complement of the vertex.
\end{Remark}

\begin{proof}
Let $Y$ be a component of the fixed point set with $\Phi (Y) = 0$.  We
use the the notation of Lemma~\ref{lemma_model} and
Remark~\ref{split}.

By Lemma~\ref{lemma_model} the zero level set $\Phi \inv (0)$ near $Y$
is isomorphic to
$$
\left\{ [q, (\vec{z}_1, \ldots, \vec{z}_k)] \left| \sum _{i=1}^{s}
\kappa _i |\vec{z}_i|^2 = \sum _{i=s+1}^{k} \kappa _i |\vec{z}_i|^2
\right. \right\} \simeq P\times _G \overset \circ c(S^+ \times S^-),
$$
where $S^+ = \{z \in V^+ \mid \sum _{i=1}^{s} \kappa _i |\vec{z}_i|^2 =
1\}$, and $S^-$ is defined similarly.  
Therefore the reduced space $\Phi \inv (0) /S^1$ near the stratum $Y$
is $P\times _G \overset \circ c (S^+ \times_{S^1} S^-)$ where the
action of $S^1$ on $S^+ \times S^- \subset V^+ \times V^-$ is defined
by the weights $\kappa_1, \ldots \kappa _k$.   
\end{proof}

\section{The perturbed quotient}

In this section we will construct an orbifold $\tred$,
which we call the {\bf perturbed quotient}, together with a map $f:
\tred \to \red$.  The perturbed quotient has two key properties: it is
straightforward to explicitly compute its cohomology ring; and $f$
induces a pairing preserving isomorphism between the cohomology ring
of the perturbed quotient and the intersection cohomology (middle
perversity) of the reduced space. 
While we do explain what $f$ looks like locally  in
this section, we defer showing that $f$ induces an isomorphism 
in cohomology to the last section.

\subsection*{The key idea}

The key idea that makes this work is an observation due to Yi Hu (this
observation was made in the context of algebraic actions on
projective varieties \cite{Hu}):  If $0$ is a singular value of an
$S^1$ moment map $\Phi: M\to \R$ and $0$ lies in the interior of the
image $\Phi (M)$ then for each component $Y$ of the fixed point set
$M^{S^1}$ with $\Phi (Y) = 0$ there exists a regular value $\epsilon
\in \R$ of $\Phi$ and a neighborhood $U$ of $Y$ in $M$ so that there is a
natural isomorphism
$$IH_{\bm} ^* \left((\Phi
\inv (0) \cap U)/S^1\right) \simeq H^* \left((\Phi \inv (\epsilon)
\cap U)/S^1\right).$$

\begin{Lemma}
We use the notation of Lemma~\ref{lemma_model}, Remark~\ref{split} and
Proposition~\ref{prop_simple}.  Fix a component $Y$ of the fixed point
set. Consider the associated bundle $P\times _G (V^+\times V^-)$,
together with the moment map
$$
 \mu ([p, z^+, z^-]) = |z^+|^2  - |z^-|^2,
$$
where $z^+ = \sum _{i=1}^{s}\vec{z}_i\in V^+$,$z^- = \sum
_{i=s+1}^{k}\vec{z}_i\in V^-$, $|z^+|^2 = \sum _{i=1}^{s}\kappa_i
|\vec{z}_i|^2$, and $|z^-|^2 = \sum _{i=s+1}^{k}\kappa_i
|\vec{z}_i|^2$.   

For any $\epsilon > 0$, $\pm\epsilon$ are regular values of $\mu$ and
$$
\mu \inv (\pm \epsilon)/S^1 = P\times _G X^\pm
$$
where 
$$
X^+ \simeq S^+ \times _{S^1} V^-, \quad\text{a $V^-$ bundle over the
weighted projective space $S^+/S^1$}.
$$
$$
X^- \simeq V^+ \times _{S^1} S^-, \quad\text{a $V^+$ bundle over the
weighted projective space $S^-/S^1$}.
$$
Here as in Proposition~\ref{prop_simple}, $S^\pm$ is the unit sphere
in $V^\pm$, $S^\pm := \{ z\in V^\pm \mid |z^\pm|^2 = 1\}$, and the
$S^1$ action on $V^+\oplus V^-$ is given by the weights $\kappa_1 ,
\ldots, \kappa_k$.
\end{Lemma}

\begin{proof}
Note first that the only fixed point in $V^+\times V^-$ under the
$S^1$ action is $(0,0)$, and that $\mu (0,0) = 0$.  Therefore any
$\alpha \not = 0$ is a regular value of $\mu$.  Next assume that $Y$
is a point; in this case  $P\times _G (V^+\times V^-)$ is simply $V^+\times
V^-$.  Then $\mu \inv (\epsilon) = \left\{ (z^+, z^-)\mid |z^+|^2 -
|z^-|^2 = \epsilon \right\} = \left\{ (z^+, z^-)\mid |z^+|^2 = |z^-|^2
+ \epsilon \right\}$.  If $\epsilon >0$, we have an $S^1$-equivariant
diffeomorphism $S^+\times V^- \to \mu \inv (\epsilon)$ given by
$(\zeta, w)\mapsto (\sqrt{\epsilon + |w|^2} \zeta, w)$.  Therefore
$\mu \inv (\epsilon)/S^1 = S^+ \times _{S^1} V^-$ for $\epsilon >0$.
Since the norm $|(z^+, z^-)|^2 = |z^+|^2 + |z^-|^2$ is $G\times
S^1$-invariant by construction (see Lemma~\ref{lemma_model} and
Remark~\ref{split}) the claim follows.
\end{proof}

Observe that $S^+\times _{S^1} S^-$ is the sphere bundle of the vector
bundle $X^+ = S^+\times _{S^1} V^-$ over the weighted projective space
$S^+/S^1$.  Therefore, by collapsing the zero section of the bundle
$X^+\to S^+/S^1$ to a point we obtain a map from $X^+$ to the cone
$\overset \circ c (S^+ \times_{S^1} S^-)$, which is a diffeomorphism
off the zero section onto the cone minus the vertex. 

We will see later on  that  if $\dim V^+ \leq \dim V^-$, then
$$
IH_{\bm}^* (\mu \inv (0)/S^1) = IH_{\bm}^* (P\times _G \overset \circ
c (S^+ \times_{S^1} S^-))  \simeq H^* (P\times _G X^+) = H^* (\mu \inv
(\epsilon )/S^1)
$$
for any $\epsilon >0$.  Similarly, if $\dim V^+ \geq \dim V^-$ 
then
$$
IH_{\bm}^* (\mu \inv (0)/S^1) \simeq H^* (\mu \inv (\epsilon
)/S^1)\quad \text{for any }\quad \epsilon < 0.
$$
In the algebraic category, these isomorphisms follow from the fact
that the collapsing map is a {\bf small resolution}, and a fact that
small resolutions induce isomorphisms in cohomology.  This isomorphism
is also valid in the symplectic context, as we prove in the next
section.

Note  that $\dim V^+ \leq \dim V^-$ if and only if the index of $Y$
as a critical manifold of the Bott-Morse function $\Phi$ is 
at most $\frac{1}{2} (\dim M - \dim Y)$.  Unfortunately we
cannot expect such inequalities to hold globally, that is, if $0$ is a
singular value of the moment map $\Phi: M\to \R$ we {\bf should not}
expect $IH_{\bm}^* (\Phi \inv (0)/S^1) = H^* (\Phi \inv (\epsilon
)/S^1)$ for some $\epsilon \not = 0$: it may well happen that at one
component $Y$ we would need to shift the value of the moment map
down and at another component to shift the value up in order to obtain
a resolution of the singularities of the reduced space at zero.

\subsection{The Construction}

Let the circle $S^1$ act effectively on a compact connected symplectic
manifold $M$ with a moment map $\Phi: M \to \R$ so that $0$ is in the
interior of the image $\Phi (M)$.
We will now construct a
Morse-Bott function $\tilde\Phi: M \to \R$ and an $S^1$ equivariant
map $f: \tPhi\inv (0)/{S^1} \to \Phi\inv (0)/{S^1}$ with the following
properties.
\begin{itemize}
\item The  critical points of $\tPhi$ are exactly the fixed points
of $S^1$ on $M$. 
\item $0$ is a regular value of $\tPhi$. 
\item The map $f:\tPhi \inv (0)/S^1 \to \Phi \inv (0)/S^1$
induces an isomorphism in cohomology  $IH^*_\bm(\Phi\inv(0)/S^1)
\cong H^*(\tPhi\inv(0)/S^1)$.
\end{itemize}

\begin{definition}
We call the subquotient $\tred := \tPhi \inv(0)/S^1$ the {\bf
perturbed quotient}.
\end{definition}
The first two properties guarantee that $\tred$ is an orbifold, and
that it is possible to compute the cohomology ring $H^*(\tred)$ in a
fairly straightforward manner.  This will be treated explicitly in
the next subsection.  As we mentioned earlier, the last property will
not be proved in this section.  However, we will prove
Lemma~\ref{lemma_smallres}, which we will later see is sufficient to
construct this isomorphism.

For each critical manifold $Y_i$ of $\Phi$ in $\Phi \inv (0)$,
there is a neighborhood $U_i$ of $Y_i$ in $M$ 
which is equivariantly isomorphic to the model 
$$
P_i \times_{G_i} (V_i^+\times V_i^-)
$$  
where the principal bundle $G_i \to P_i \to Y_i$ and the vector spaces
$V_i^+$, $V_i^-$ are as in the preceding section.  We may assume that
the $U_i$'s  for distinct critical manifolds do not intersect.
There exists $\delta >0$ so  that 
$0$  is the only critical value of $\Phi$ in $(-\delta, \delta)$
and $U_i$ is the image of the set 
$$
P_i \times_{G_i} (\{ (z_i^+, z_i^-)\mid |z_i^+|^2 + |z_i^-|^2< 3\delta  \}.
$$

Therefore, we will simply give our construction on the vector space
$V^+\times V^-$. 
As long as our definition of $\tPhi$ and $f$  are $G$-invariant,
these construction can be naturally extended to the local model.
Additionally, as long as   $\tPhi = \Phi$ and $f$ is the identity 
outside the set
$\{ (z_i^+, z_i^-)\mid |z_i^+|^2 + |z_i^-|^2< 3\delta \}$,
they can be extended globally by taking
$\tPhi = \Phi$ and $f = \id$ on $M \smallsetminus \cup U_i$.

Choose a smooth function $\rho: \R \to \R$ such that $\rho (t) = 1$
for all $t<\delta$, $\rho (t) = 0$ for all $t>2\delta$ and $\rho
'(t)\leq 0$ for all $t$.  Let $C = \sup |\rho ' (t)|$, and choose
$\epsilon \in \R$ so that $\epsilon \neq 0$, 
$|\epsilon | < C\inv$ and  $|\epsilon | <\delta$.  
Moreover, choose $\epsilon$ so that $\epsilon > 0$ if and only  if 
$\dim V^+ \leq \dim V^-$.
We now define our new function $\tPhi$ :
$$
\tPhi (z^+, z^-) := 
\Phi (z^+, z^-) + \epsilon \rho (|(z^+, z^-)|^2) =
|z^+|^2 - |z^-|^2 + \epsilon \rho (|(z^+, z^-)|^2).
$$ 

The norm 
$$
|(z^+, z^-)|^2 = |z^+|^2 + |z^-|^2
$$
is $G\times S^1$-invariant by construction (see
Lemma~\ref{lemma_model}, Proposition~\ref{prop_simple} and the
subsequent discussion).  Therefore the function $\rho (|(z^+,
z^-)|^2)$, and hence also the function $\tPhi$, is $G\times S^1$
invariant.  Moreover, for $(z^+, z^-)$ with $|(z^+,z^-)|^2 >2\delta$,
$ \tPhi(z^+, z^-)= |z^+|^2 - |z^-|^2 .$
Therefore
$$
\tPhi \inv (0) \cap \{|(z^+,z^-)|^2 > 2\delta \} = \Phi \inv (0)
\cap \{|(z^+,z^-)|^2 > 2\delta \}.
$$
In contrast, for $(z^+, z^-)$ with $|(z^+,z^-)|^2 < \delta$, 
$ \tPhi(z^+, z^-)= |z^+|^2 - |z^-|^2 + \epsilon. $
Thus $(0,0)$ is a nondegenerate critical point, and
$$
\tPhi \inv (0) \cap \{|(z^+,z^-)|^2 < \delta \} = \Phi \inv (- \epsilon)
\cap \{|(z^+,z^-)|^2 < \delta\}.
$$
(Note that $|\epsilon| < \delta$ guarantees that $\Phi \inv (-
\epsilon) \cap \{|(z^+,z^-)|^2 < \delta\} \not = \emptyset$.)
Moreover, $(0,0)$ is the only critical point of $\tPhi$, because
\begin{equation*}
\begin{split}
d \tPhi & = d |z^+|^2 - d |z^-|^2 + \epsilon d\rho (|(z^+,
z^-)|^2) \\
& = \left(1 + \epsilon \rho' (|(z^+,z^-)|^2)\right)d |z^+|^2  -
  \left(1 - \epsilon \rho' (|(z^+,z^-)|^2)\right)d |z^-|^2  
\end{split}
\end{equation*}
and $|1 \pm \epsilon \rho'| \geq 1 - |\epsilon| (\sup |\rho ' (t)|) >
0$, since $|\epsilon| (\sup |\rho ' (t)|) < 1$ by the choice of
$\epsilon$.  It follows that $\tPhi$ is a Bott-Morse function, and
that $0$ is a regular value of $\tPhi$ (since $\tPhi (0,0) = \epsilon
\not = 0$).

\begin{definition}
Let $X= X^r\sqcup \coprod Y_i$ be a simple stratified space.  A {\bf
resolution} $h:\tilde X \to X$ is a continuous surjective map from a
smooth orbifold $\tilde X $ such that $h\inv (X^r)$ is dense in
$\tilde X$ and $h: f\inv (X^r)\to X^r$ is a diffeomorphism.
\end{definition}

We will now construct a resolution $f :  \tred \to \red$.
We start by considering a $G\times S^1$-equivariant
map $\psi : V^+\times V^- \to \Phi \inv (0)$ defined by
\begin{equation}\label{eq_fiber}
\begin{split}
\psi (z^+, z^-) &= ( \left(\frac{|z^-|^2}{|z^+|^2}\right)^{1/4} z^+,
\left(\frac{|z^+|^2}{|z^-|^2}\right)^{1/4} z^-) \quad\text{if} \quad
z^+, z^- \not = 0\\
\psi (0, z^-) &= \psi (z^+,0) = (0, 0).\\
\end{split} 
\end{equation} 
We let $f: \tred \to \red$ be $G$-equivariant map induced
by the restriction  $\psi|_{\tPhi\inv(0)} : \tPhi\inv (0) \to \Phi \inv (0)$.

To prove that $f$ is a resolution, it is enough to show that $\psi |_{\tPhi\inv
(0) \smallsetminus \psi\inv (0,0)} : \tPhi\inv (0) \smallsetminus
\psi\inv (0,0) \to \Phi \inv (0) \smallsetminus \{(0,0)\}$ is a
diffeomorphism.  It follows from~(\ref{eq_fiber}) that for $(0,0) \not
= (z^+, z^-) \in \Phi\inv (0)$,
\begin{equation} \label{eq_fiber_psi}
\psi \inv (z^+, z^-) = \{ (\lambda z^+, \lambda^{-1} z^-) \mid
\lambda >0 \}.
\end{equation}
Consequently $\psi|_{\tPhi\inv (0) \smallsetminus \psi\inv (0,0)} :
\tPhi\inv (0) \smallsetminus \psi\inv (0,0) \to \Phi \inv (0)
\smallsetminus \{(0,0)\}$ is one-to-one and onto.  Therefore it
remains to prove that $d \psi |_{T (\tPhi\inv (0) \smallsetminus
\psi\inv (0,0)}$ is one-to-one, or, equivalently, that for any $(z^+,
z^-) \in \tPhi\inv (0) \smallsetminus \psi\inv (0,0)$ 
$$
0 = \ker d\psi \cap T_{(z^+, z^-)} \tPhi\inv (0) = \ker d\psi \cap
\ker d\tPhi.
$$
By (\ref{eq_fiber_psi}), the kernel of $d\psi$ at $(z^+, z^-)$ is
spanned by the vector $\left.\frac{d}{d\lambda}\right|_{\lambda = 1}
(\lambda z^+,\lambda^{-1} z^-)$.  Thus  it remains to show that for any 
$(z^+, z^-) \in \tPhi \inv (0) \smallsetminus \psi \inv (0,0)$ 
we have 
$$
\left. \frac{d}{d\lambda}\right|_{\lambda = 1} \tPhi (\lambda
z^+,\lambda^{-1} z^-) \not = 0.
$$
Now
\begin{multline*}
\left. \frac{d}{d\lambda}\right|_{\lambda = 1} \left( |\lambda z^+|^2
- |\lambda \inv z^-|^2 + \epsilon \rho (|(\lambda z^+|^2+ |\lambda\inv
z^-|^2)\right) = \\ 
\left. \left( 2 \lambda |z^+|^2 + 2 \lambda ^{-3} |z^-|^2 + \epsilon
\rho' (|(\lambda z^+|^2+ |\lambda\inv z^-|^2)(2 \lambda |z^+|^2 - 2
\lambda ^{-3} |z^-|^2) \right) \right| _{\lambda =1} = 
\\ 
2 \left(|z^+|^2 + |z^-|^2\right) \ + \  
2 \epsilon \rho'(|z^+|^2 + |z^-|^2) (|z^+|^2 - |z^-|^2).
\end{multline*}
For $(z^+,z^-) \in \tPhi\inv(0)$,  we have
$|z^+|^2 - |z^-|^2 = -\epsilon \rho(|(z^+,z^-)|) $. 
Since $ \rho'(t) \leq 0$ for all $t$, $-\epsilon ^2 \rho (t) \rho'
(t)\geq 0$ for all $t$. Moreover, $(z^+,z^-) \neq (0,0)$. Hence 
\begin{multline*}
\left. \frac{d}{d\lambda}\right|_{\lambda = 1} \tPhi (\lambda
z^+,\lambda^{-1} z^-) =  \\
2 \left(|z^+|^2 + |z^-|^2\right)
 - 2 \epsilon \rho'(|z^+|^2 + |z^-|^2) \epsilon \rho (|z^+|^2 + |z^-|^2) 
\ \ \geq  \ \ 
2 \left(|z^+|^2 + |z^-|^2 \right)
 > 0.
\end{multline*}
Thus, we have proved the following.

\begin{Lemma} \label{lemma_smallres}
Let a circle $S^1$ act on a symplectic manifold $M$ with a moment
map $\Phi: M \to \R$ so that $0$ is in the interior of the
image  $\Phi(M)$.  
Let $\tPhi: M \to \R$ and 
$f :\tred  = \tPhi \inv (0)/S^1\to \red = \Phi \inv (0)/S^1$
be constructed as above.

Then $f$ is a resolution.
Moreover, for each singular stratum $Y$ of $\red$ there exist:
\begin{itemize}
\item an even dimensional orbifold vector bundle $E\to N$ over a
compact orbifold $N$ with a sphere bundle $L\to N$ such that
$$
  \dim N \leq \frac{1}{2} \dim E - 1,
$$
\item a principal $G$ bundle $P \to  Y$,
\item an action of $G$ on $E$ by vector bundle maps
\item an isomorphism from a neighborhood of the vertex section of the
cone bundle $P\times _G \conec (L) \to Y$ to a neighborhood $U$ of $Y$
in $\red$,
\item an isomorphism from a neighborhood of the zero section of the
vector bundle $P\times _G E \to P\times _G N$ to the neighborhood
$f\inv (U)$ of $f\inv (Y)$ in $\tred$ such that the diagram
$$
\begin{CD}
P\times _G E @<<<  f\inv (U) @>>> \tred\\
@V{h}VV                     @VV{f}V                 @VV{f}V\\
P\times _G \conec (L) @<<< U @>>> \red
\end{CD}
$$
commutes.  Here the map $h$ is induced by the natural blow-down map
$E\to \conec (L)$ taking the zero section to the vertex.
\end{itemize}
\end{Lemma}

Notice that there is no reason to suspect that the perturbed
quotient possesses a symplectic structure.

\begin{Remark}
Morally, the above Lemma should be read as a claim that
$f:\tred \to \red$ is a small resolution (cf. \S 6.2 of \cite{GM2}).
\end{Remark}

\subsection{Computation of the  cohomology of the perturbed quotient}

We can now compute the cohomology of the perturbed quotient by
adapting techniques used to compute the cohomology ring of a
symplectic quotient at a regular value.

We begin by reviewing those techniques.  Let a circle $S^1$ act on a
compact connected symplectic manifold $M$ with a moment map $\Phi$.
Assume that $0$ is a regular value.  There is a natural restriction
from $H_{S^1}^*(M; \R)$, the equivariant cohomology of $M$, to
$H_{S^1}^*(\Phi\inv(0); \R))$, the equivariant cohomology of the
preimage of $0$.  Since $0$ is a regular value, the stabilizer of
every point in $\Phi\inv(0)$ is discrete.  Therefore, there is a
natural isomorphism from $H^*_{S^1}(\Phi\inv(0),\R)$ to the
$H^*(\red)$, the ordinary cohomology of the symplectic quotient $\red
: \Phi\inv(0)/S^1$.  The composition of these two maps gives a natural
map, $\kappa : H_{S^1}^*(M) \to H^*(\red)$, called the {\bf Kirwan
map}.

\begin{Theorem}(Kirwan, [K]) 
Let a circle $S^1$ act on a compact connected symplectic manifold $M$
with a moment map $\Phi$ so that $0$ is a regular value.  The Kirwan
map $\kappa : H_{S^1}^*(M; \R) \to H^*(\red ; \R)$ is surjective.
\end{Theorem}

Thus, assuming we know the ring structure on $M$, the ring structure
on $\red$ can be computed from the kernel of $\kappa$.  By Poincare
duality, to compute the kernel it is enough to compute the integral of
$\kappa(\alpha)$ over the reduced spaces for every equivariant
cohomology class $\alpha$ on $M$.  We take one formula for this
integral from Kalkman \cite{Ka}; slightly different but morally
equivalant formulas were proved by Wu \cite{Wu} and a more general
version by Jeffrey-Kirwan \cite{JK}. See also \cite{Ka-Gu}.  All of
these results were inspired by a paper of Witten  \cite{Wi}.

\begin{Theorem}
Let a circle $S^1$ act on a compact connected symplectic manifold $M$
with a moment map $\Phi$ so that $0$ is a regular value. Let $\cF^+$
denote the set of components $F$ of the fixed point set $M^{S^1}$ such
that $\phi(F) > 0$.

Given an equivariant cohomology class $\alpha \in H^*_{S^1}(M)$, 
the integral of $\kappa(\alpha)$ over $\red$ is
given by the formula
$$
\int_{\red} \kappa(\alpha) = \Res _0 \sum_{F \in \cF^+} \int_F \frac{
i_F^*(\alpha) }{e_F} ,
$$
where $e_F$ denotes the equivariant Euler class of the normal bundle
of $F$.
\end{Theorem}

The right hand side of this formula requires some explanation.  The
map $i_F^*$ is simply the restriction to $F$.  The equivariant
cohomology ring $H^*_{S^1}(F)$, is naturally isomorphic to $H^*(F)[t]$.
The equivariant Euler class $e_F$ is invertable in
the localized ring $H^*(F)(t)$; thus,
$\frac{i_F^*(\alpha)}{e_F}$ is an element of this ring.  The integral
$\int_F : H^*(F)(t) \to \R(t)$ acts by integrating each
coefficient in the series.   Finally, $\Res _0$ denotes
the operator which returns the coefficient  of $t\inv$.

An alternative way of computing the kernel is given by a theorem of
Tolman and Weitsman.
\begin{Theorem}[Tolman-Weitsman, \cite{TW}] 
Let the circle $S^1$ act on a compact connected symplectic 
manifold $M$ with a moment map $\Phi: M \to \R$
so that $0$ is a regular value.

Let $\cF^+$ denote the set of components $F$ of the fixed point set
such that $\Phi(F) > 0$; let $\cF^-$ denote the set of components $F$
of the fixed point set such that $\Phi(F) < 0$.  Define
$$
K_\pm := \{ \alpha \in H_{S^1}^*(M) \mid \alpha|_{F} = 0 \ \ \forall \ \ 
F \in {\cF}^\pm\}.$$

The kernel of the Kirwan map is $K_+ \oplus K_-$.
\end{Theorem}

In our case, closely analogous propositions are true.

\begin{Proposition}\label{compute}
Let the circle $S^1$ act on a compact connected symplectic manifold
$2n$ dimensional manifold $M$ with a moment map $\Phi: M \to
\R$. Assume 
that $0$ is in the interior of $\Phi(M)$.  Let $\tred$ denote the
perturbed quotient.

Then there is a surjective  ring homomorphism
$\kappa: H^*_{S^1}(M) \to H^*(\tred)$.
Moreover, 
\begin{itemize}
\item
The kernel of $\kappa$ is $K_+ \oplus K_-$, where
$K_\pm := \{ \alpha \in H_{S^1}^*(M) \mid \alpha|_{F} = 0 \ \forall \ \ 
F \in {\cF}^\pm\}. $
\item
Given an equivariant cohomology class  
$\alpha \in  H^*_{S^1}(M)$, 
the integral of $\kappa(\alpha)$ over $\red$ is given by the formula
$$\int_{\red} \kappa(\alpha) = 
\Res_0 \sum_{F \in \cF^+} \int_F \frac{ i_F^*(\alpha) }{e_F} ,$$
where  $e_F$ denotes the equivariant Euler class of the 
normal bundle of $F$.

Here, et $\cF^+$ denotes  the set of components  $F$
of the fixed point set $M^{S^1}$  such that either
\begin{enumerate}
\item $\Phi(F) > 0$ or
\item  $\Phi(F) = 0 $ and $ 2 \index  F + \dim F \leq \dim M$.
\end{enumerate}
Additionally,  $\cF^- $  denotes all other components of the fixed
point set, $\cF^- = \cF \setminus \cF^+$.
\end{itemize}
\end{Proposition}

The reason that this proposition is true is that the perturbed quotient
$\tred$ is defined in a way very similar to the ordinary reduced
space.

Thus, for example, Kalkman's formula follows immediately from the fact
that there exists a smooth invariant function $\Phi : M \to \R$ so
that $0$ is regular and $\tred$ is defined to be $\tPhi\inv(0)/S^1$.
His proof relies only on the fact that $\tPhi\inv([0,\infty))$ is a
manifold with boundary.  Thus, one only need note that $\cF^+$ does
indeed correspond to the components $F$ of the fixed point set such
that $\tPhi(F) > 0$.

To see that Kirwan's surjectivity holds, and the Tolman-Weitsman
formula for the kernel of $\kappa$, we must also use the fact that
$\tPhi$ is a Morse-Bott function and that its critical points are
exactly  the fixed points of the action.  This is sufficient to
prove both theorems, as was pointed out in \cite{TW}.

\section{The isomorphism}

The goal of this section is to prove that the intersection cohomology
of the symplectic quotient by a Hamiltonian circle actoin is
isomorphic to the (ordinary) cohomology of the perturbed quotient.
Because we have already computed the cohomology of the perturbed
quotient in Proposition \ref{compute}, this will allow to obtain the
description of the intersection cohomology of the symplectic quotient,
and thus prove our main theorems.  More precisely, we will be done
once we have proved the following.

\begin{Theorem}\label{theorem-iso}
Let the circle $S^1$ act on a compact connected symplectic manifold
$M$ with moment map $\Phi: M \to \R$ so that $0$ is in the interior of
$\Phi(M)$.  Let $\red := \Phi\inv(0)/S^1$ denote the reduced space and
let $\tred$ denote the perturbed reduced space.  There is a natural
pairing preserving isomorphism between the intersection cohomology of
the symplectic quotient $\red$ and the cohomology if the perturbed
quotient $\tred$.

More precisely there exists an isomorphism $\psi: H^* (\tred) \to IH^*
(\red)$ of graded vector spaces such for any $\alpha \in H^p (\tred)$
and $\beta \in H^q (\tred)$ with $p+q = \dim \tred$  we have 
$$
 \int _{\tred} \alpha \cup \beta = \int _{\red} \langle \psi (\alpha),
 \psi (\beta) \rangle .
$$

\end{Theorem}

Instead of trying to construct the isomorphism between the
intersection cohomology of the reduced space and the (ordinary)
cohomology of the perturbed quotient directly, we will introduce a new
complex $A^*_\bm(\tred) = A^*_\bm(f:\tred \to \red) $ and show that
the cohomology of $A^*_\bm (\tred \to \red) $ is naturally isomorphic
to both $H^* (\tred )$ and $IH_\bm^* (\red)$.

\begin{definition}
 Let $f: \tX \to X$ be a resolution of a simple stratified space.  Let
$X^r$ be the top stratum of $X$, $\tX^r$ be its preimage $f\inv(X^r$),
and let $\iota: \tilde{X}^r \hookrightarrow \tilde{X}$ denote the
inclusion.  By construction, there are maps of complexes $f^* :
I\Omega _{\bm}^* (X) \to \Omega^* (\tilde{X}^r )$ and $\iota ^* :
\Omega ^*(\tilde{X}) \to \Omega^* (\tilde{X}^r )$.  We define the
complex of {\bf resolution forms}
\begin{equation}\label{the_complex}
A_\bm^\bullet (\tX) = A_\bm^\bullet (f: \tilde{X} \to X) := f^*
\left( I\Omega _{\bm}^\bullet (X)\right ) \cap \iota^* \left(\Omega
^\bullet (\tilde{X}) \right)
\end{equation}
\end{definition}
Note that $f^*$ and $\iota^*$ are both injective.  Therefore we may
think of a resolution form as an intersection form on $X^r$ which
extends to a  globally defined form on $\tX$.  This gives us the
inclusions of complexes $A_\bm^\bullet (\tX) \to I\Omega _{\bm} ^\bullet (X)$
and $ A_\bm^\bullet (\tX) \to \Omega^\bullet(\tilde{X})$, which induce
maps in cohomology $j:H^* (A_\bm^\bullet (\tX)) \to IH^*_\bm (X)$ and
$i:H^* (A_\bm^\bullet (\tX)) \to H^* (\tX)$.  Note that the graded
vector space $H^*(A_\bm (\tX)$ has a pairing defined by taking the
exterior product of the representatives of the classes and then
integrating the product over $\tX$.  Clearly the maps $i$ and $j$ are
pairing preserving.  Thus, to prove Theorem~\ref{theorem-iso}  it
is enough to show that the maps $i$ and $j$ are isomorphisms.

\subsection{Local issues}

We start with a simple calculation.
\begin{Lemma}\label{lemma_calculation}
Let $E\to N$ be an even dimensional orbifold vector bundle over an
orbifold $N$, and let $L$ denote the sphere bundle of $E$.  Then the
obvious blow-down map $f: E \to \conec (L)$ is a resolution.  If $\dim
N \leq \frac{1}{2} E -1$ then the maps
$$
  A_\bm^\bullet (E) \hookrightarrow \Omega^\bullet (E)
$$
and 
$$
A_\bm^\bullet (E) \hookrightarrow I\Omega^\bullet_\bm (\conec (L))
$$
induce isomorphisms in cohomology.
\end{Lemma}
To prove the Lemma we will need the following technical observation.

\begin{Lemma}\label{prop_cyl}
Let $L$ be an orbifold.  Let $\alpha$ be a closed $k$ form on the
cylinder $L\times (0,\infty)$ which vanishes on $L\times (0, a)$ for
some $a$.  Then $\alpha = d\beta$ for some $k-1$ form $\beta$ which
also vanishes on  $L\times (0, a)$.
\end{Lemma}
\begin{proof}
Consider first the case where $L$ is a point and $\alpha = f(r) \,dr$
is a 1 form on $(0, \infty)$.  Then $\beta = \int _0 ^r f (s)\, ds$.

In general, if $L$ is not a point, the $k$-form $\alpha$ has to be of
the form $f(r)\wedge dr$ where $f(r)$ is in $\Omega^{k-1} (L)$ for
each $r\in (0,\infty)$.  Let $\beta = \int _0 ^r f(s) \, ds$.
\end{proof}

\begin{proof}[Proof of Lemma~\ref{lemma_calculation}]
Note first that the middle perversity of the vertex $*$ of the cone
$\conec (L)$ is $\bm (*) = \frac{1}{2} \dim E -1$, since $E$ is even
dimensional.   

Recall that $\conec (L)$ is a stratified space with two strata: the
vertex $*$ and the complement $ L\times (0, \infty)\simeq
E\smallsetminus N$.  We may choose the tubular neighborhood $T$ to be
any neighborhood of the vertex $*$ of the form $L\times (0,a)/\sim$
for some $a$.

It follows from the definitions that
$$
I\Omega^q_\bm (\conec (L)) = 
\begin{cases}
\Omega ^q (E\smallsetminus N)&
\text{for $q< \bm (*)$}\\ 
\left\{ \alpha \in \Omega ^q (E\smallsetminus N) \mid d\alpha |_T =
0\right\} & \text{ for $q = \bm (*)$}\\ 
\left\{ \alpha \in \Omega ^q (E\smallsetminus N) \mid\alpha |_T = 0
\text{ and } d\alpha |_T = 0\right\} & \text{ for $q > \bm (*)$}
\end{cases}
$$
Consequently
$$
A^q_\bm (E) = 
\begin{cases}
\Omega^q (E)&
\text{for $q< \bm (*)$}\\ 
\left\{ \alpha \in \Omega ^q (E) \mid d\alpha |_T =
0\right\} & \text{ for $q = \bm (*)$}\\ 
\left\{ \alpha \in \Omega ^q (E) \mid\alpha |_T = 0
\text{ and } d\alpha |_T = 0\right\} & \text{ for $q > \bm (*)$}
\end{cases}
$$
The map $A^q_\bm (E) \to I\Omega_\bm^q (\conec (L))$ is induced by the
restriction from $\Omega ^* (E) $ to $\Omega ^* (E\smallsetminus N)$.

It follows from Lemma~\ref{prop_cyl} that for $q\leq \bm (*)$
the map $H^q (A^\bullet _\bm (E)) \to H^q (E)$ is an isomorphism and
that $H^q (A^\bullet _\bm (E)) =0$ for $q> \bm (*)$.  Since $H^* (E) =
H^* (N)$ and since $\bm (*) \geq \dim N$ by assumption, the map $H^q
(A^\bullet _\bm (E)) \to H^q (E)$ is an isomorphism for all $q$.

Similarly, 
$$
IH^q_\bm (\conec (L)) = 
\begin{cases}
H^q (E\smallsetminus N)&
\text{for $q\leq \bm (*)$}\\ 
0  & \text{ for $q >\bm (*)$}
\end{cases}.
$$
Consider the Gysin sequence
$$
\cdots \to H^{q- \lambda} (E) \to H^q (E) \to H^q (E\smallsetminus N)
\to H^{q -\lambda +1} (E) \to \cdots
$$
where $\lambda = \dim E - \dim N$.  Since for $q -\lambda + 1 \leq -1$
(i.e., for $q \leq \lambda -2$) we have $ H^{q -\lambda +1} (E) = 0 =
H^{q -\lambda } (E) $, the pull-back map $H^q (E) \to H^q
(E\smallsetminus N)$ is an isomorphism.  In particular the pull-back
is an isomorphism for for $q \leq \bm (*)= \frac{1}{2} \dim E - 1 =
\dim E - (\frac{1}{2} \dim E - 1) - 2 \leq \dim E - \dim N - 2 =
\lambda - 2$.
\end{proof}

\begin{Proposition}
\label{local-iso}
Let the circle $S^1$ act on a compact connected symplectic manifold
$M$ with moment map $\Phi: M \to \R$. Assume that $0$ is in the interior of
the image $\Phi(M)$.  Let $\red := \Phi\inv(0)/S^1$ denote the reduced
space and let $f: \tred \to \red$ denote its resolution by the
perturbed quotient.  

There exists a cover $\cU$ of $\red$ such that the natural inclusions
$A_\bm^* (f\inv(U_{\alpha_1}) \cap \cdots \cap f\inv(U_{\alpha_k}) )
\to I\Omega ^*_{\bm} (U_{\alpha_1} \cap \cdots U_{\alpha_k})$ and $
A_\bm^* (f\inv (U_{\alpha_1}) \cap \cdots \cap f\inv(U_{\alpha_k}))
\to \Omega^*( f\inv(U_{\alpha_1}) \cap \cdots \cap f\inv(U_{\alpha
_k}))$ induce isomorphisms in cohomology for all $k$-tuples
$\{ U_{\alpha_1}, \ldots, U_{\alpha_k} \}$ of elements of $\cU$.
\end{Proposition}

\begin{proof}
We have seen in Proposition~\ref{prop_simple} that $\red = \red^r
\coprod Y_i$ where $Y_i$ are compact manifolds.  Further, for each
singular stratum $Y$ there exists a tubular neighborhood $\tT$ of
$Y$, the fiber bundle $\conec (L)\to \tT\stackrel{\pi}{\to} Y$, and
the map $r:\tT\to [0,1)$ (c.f.\ Remark~\ref{rmrk_t}).  Recall also
that by Lemma~\ref{lemma_smallres} we may assume that $\tT = P\times
_G \conec (L)$, that $f\inv (\tT) = P\times _G E$ and that $f: P\times
_G E \to P\times _G \conec (L)$  is induced by the obvious blow-down
map $E\to \conec (L)$.

We take $U_0 = \red \smallsetminus \cup r_i \inv ([0, 1/2]) = \red
\smallsetminus \overline{T}_i$.  Since each singular stratum $Y$ is a
compact manifold, it possesses a finite good cover $\{V_\alpha \}$.
Moreover we may assume that $\pi \inv (V_\alpha ) \simeq \conec
(L)\times V_\alpha $.  We take $U_\alpha := \pi \inv (V_\alpha )
\cap \tT \subset \tT$.  This give us a cover $\cU$ of $\red$.

Note that by construction for a $k$-tuple $\{ U_{\alpha_1}, \ldots,
U_{\alpha_k} \}$ of elements of $\cU$ we either have that
$U_{\alpha_1} \cap \cdots \cap U_{\alpha_k}$ does not intersect any
singular stratum $Y$ (in which case $f\inv (U_{\alpha_1}) \cap \cdots
\cap f\inv(U_{\alpha_k})$ and $U_{\alpha_1} \cap \cdots \cap U_{\alpha_k}$
are diffeomorphic) or there is a unique stratum $Y$ such that $Y \cap
U_{\alpha_1} \cap \cdots \cap U_{\alpha_k} \not = \emptyset$.  In the
latter case $U_{\alpha_1} \cap \cdots \cap U_{\alpha_k} \simeq D \times
\conec (L)$, $f\inv (U_{\alpha_1}) \cap \cdots \cap f\inv(U_{\alpha_k})
\simeq D\times E$ and $f: f\inv (U_{\alpha_1}) \cap \cdots \cap
f\inv(U_{\alpha_k}) \to U_{\alpha_1} \cap \cdots U_{\alpha_k}$ is
equivalent to the map $ h \times id: E\times D \to \conec (L)\times D$,
where $D$ is a disk in $Y$ and $h: E \to \conec (L)$ is the resolution.

Given a disk $D$ and a set $X$ we have an inclusion $\iota : X
\hookrightarrow X \times D$, $\iota (x) = (x, 0)$.  Clearly the
diagram
$$
\begin{CD}
E @>{\iota}>> E \times D\\
@V{h}VV	                    @VV{h\times id}V\\
\conec (L)   		@>{\iota}>> \conec (L) \times D
\end{CD}
$$
commutes. Since $h: E\to \conec (L)$ and $h\times id : E \times D \to
\conec (L) \times D$ are resolutions, we have a commutative diagram of
complexes
\begin{equation}\label{eq_diagram}
\begin{CD}
\Omega^* (E) @<{\iota^*}<<\Omega^* (E \times D) 
@= \Omega^* (f\inv (U_{\alpha_1}) \cap \cdots \cap f\inv(U_{\alpha_k}))\\
@AAA @AAA\\
A^*_{\bm}(E) @<{\iota^*}<< A^*_{\bm}(E \times D)
@= A_\bm ^* (f\inv (U_{\alpha_1}) \cap \cdots \cap f\inv(U_{\alpha_k}))\\ 
@VVV @VVV\\
I\Omega_{\bm}^*(\conec (L)) @<{\iota^*}<< I\Omega_{\bm}^*(\conec
(L)\times D)
@= I\Omega^*_\bm  (U_{\alpha_1} \cap \cdots \cap U_{\alpha_k}).
\end{CD}
\end{equation}
Since the disk $D$ is contractible, the horizontal maps induce
isomorphisms in cohomology.  Since the left vertical maps induced
isomorphisms in cohomology by Lemma~\ref{lemma_calculation}, the right
vertical maps induce isomorphisms as well.
\end{proof}

\begin{Proposition}
Let the circle $S^1$ act on a compact connected symplectic manifold
$M$ with moment map $\Phi: M \to \R$. Assume that $0$ is in the interior of
the image $\Phi(M)$.  Let $\red := \Phi\inv(0)/S^1$ denote the reduced
space and let $f: \tred \to \red$ denote its resolution by the
perturbed quotient.  

The inclusions $A_\bm^* (\tred) \to I\Omega _{\bm}^* (\red)$ and $
A_\bm^* (\tred) \to \Omega^*(\tred)$ induce isomorphisms in
cohomology.
\end{Proposition}

\begin{proof}
The proof is now a standard spectral sequence argument.   

Let $\cU = \{U_\alpha \}$ be the cover of $\red$ constructed in the
proof of Proposition~\ref{local-iso}.  We now construct a continuous
partition of unity $\rho_\alpha $ subordinate to the cover $\cU$ with
the properties that
\begin{itemize}
\item the functions $\rho_\alpha $ restrict to smooth functions on
$\red^r$,
\item the functions $\rho_\alpha $ are constant along the fibers of
$\pi :T\to Y$ for all the singular strata $Y$.
\end{itemize} 
These properties ensure that 
\begin{itemize}
\item $\{f^*\rho_\alpha \}$ is a partition of unity on $\tred$
subordinate to the cover $\{f\inv (U_\alpha)\}$ of $\tred$, and that
\item for any intersection form $\gamma \in I\Omega _\bm^* (\red)$ 
or resolution form $\delta \in A_\bm^*$,
the products $\rho _\alpha \gamma$ and $f^*\rho_\alpha \delta$
are also in $I\Omega _\bm^*(\red)$ and $A^*_\bm$, respectively.
\end{itemize}
We first consider the set $U_0$ which is entirely contained in the
smooth part $\red^r$ of the quotient.  We choose $\trho _0$ to be a
smooth nonnegative function on $\red^r$ supported in $U_0$ with $\trho
_0 = 1$ in $\red \smallsetminus \cup r_i \inv ([0, 3/4))$.

By construction 
a set $U_\alpha$ with $U_\alpha \cap Y \not = \emptyset$ is of
the form $\pi \inv (V_\alpha ) \cap \tT$ where $\{V_\alpha \}$ is a
good cover of $Y$. We can choose a smooth partition of unity $\{\tau
_\alpha \}$ on $Y$ subordinate to $\{V_\alpha \}$ and also a
nonnegative smooth function $\sigma$ on $\red^r$ supported in $ \tT$
with $\sigma$ identically 1 on the $\overline{T} \cap \red^r = r\inv
([0, 1/2]) \cap \red ^r$.  Let $\trho _\alpha := \sigma (\pi^* \tau
_\alpha |_{\red})$; it extends to a continuous function on $\red$.  The
functions $\rho _\alpha := \frac{\trho _\alpha }{\sum_{\beta} \trho
_\beta }$ form the desired partition of unity.

Next we define three double complexes whose ${i,j}$'th terms are given as
follows for $j \geq 0$:
$$ 
A^{i,j}_\bm(\cU)  :=
\oplus  A^i_\bm(f\inv(U_{\alpha _0}) \cap \cdots \cap f\inv(U_{\alpha _j}))$$
$$ 
I\Omega^{i,j}_\bm(\cU) := \oplus  I\Omega_\bm^i(U_{\alpha _0}
\cap \cdots \cap U_{\alpha _j}), \quad \text{ and}
$$
$$ \tilde{\Omega}^{i,j}(\cU) := \oplus \Omega^i(f\inv(U_{\alpha _0})
\cap \cdots \cap f\inv(U_{\alpha _j})),
$$ 
where the sums are taken over all $j+1$-tuples $\{\alpha _0,
\ldots \alpha _j\}$.    
For all three complexes the differentials are given by the
de Rham and \v{C}ech differentials.

First, in order to show that the cohomology of the double complexes
are the intersection cohomology of $\red$, the cohomology of the
complex $A^i_\bm(\tred)$ and the cohomology of $\tred$ respectively,
we will consider the spectral sequences associated to the filtration
by $i$.  Since we constructed a nice partition of unity subordinate to
the locally finite cover $\cU$, for any fixed $i$ the \v{C}ech
cohomology of the sheaf $I\Omega^i_\bm$ is trivial for $j > 0$, and
for $j = 0$ consists of the global forms $I\Omega^i_\bm(\red)$.  Thus,
the spectral sequence converges at the $E_2$ term. Moreover,
$E_2^{i,j} = 0$ for $j > 0$, and $E_2^{i,0} = IH_\bm^i(\tred)$ for all
$i$.  Thus, the cohomology of the double complex
$I\Omega_\bm^{i,j}(\cU)$ is the intersection cohomology
$IH_\bm^i(\red)$.  Virtually identical arguments show that the
cohomology of the double complexes $A^{i,j}_\bm(\cU)$ and
$\tilde{\Omega}^{i,j}(\cU)$ are the cohomology of the complexes
$A^i_\bm(\tred)$ and $\Omega^i(\tred)$, respectively.

Next, in order to show that inclusions induce isomorphism from
$H^*(A_\bm(\tred))$ to $IH^*_\bm(\red)$ and from $H^*(A_\bm(\tred))$
to $H^* (\tred)$ respectively, we will consider the spectral sequences
associated with the filtration by $j$.  The double complex $
A^{i,j}_\bm(\cU) := \oplus A^i _\bm(f\inv(U_{\alpha _1}) \cap \cdots
\cap f\inv(U_{\alpha _j}))$ includes naturally into the double complex
$ I\Omega^{i,j}_\bm(\cU) := \oplus I\Omega_\bm^i(U_{\alpha _1} \cap
\cdots \cap U_{\alpha j})$.  By Proposition~\ref{local-iso} this
inclusion induces an isomorphism on the $E_1$ terms of these spectral
sequences. This implies that the inclusion induces an isomorphism on
every $E_k$.  Hence, by the proceeding paragraph, inclusion induces an
isomorphism from $H^*(A_\bm(\tred))$ to $IH^*_\bm(\red)$. An
essentially identical argument shows that the inclusion
$A^*_\bm(\tred) \to \Omega^*(\tred)$ induces an isomorphism in
cohomology.
\end{proof} 
This completes the proof of Theorem~\ref{theorem-iso}.  By combining
Theorem~\ref{theorem-iso} with Proposition~\ref{compute} we now obtain
the main result of the paper: Theorem~1 and Theorem~$1'$.

\end{document}